\documentclass[11pt,reqno]{amsart}

\usepackage{amsmath,amssymb,amsthm,mathtools}
\usepackage[margin=1.1in]{geometry}
\usepackage[T1]{fontenc}
\usepackage[utf8]{inputenc}
\usepackage{fancyvrb}
\usepackage{upquote}
\usepackage{url}
\usepackage{xcolor}
\usepackage{booktabs}
\usepackage{tikz}
\usetikzlibrary{arrows.meta}
\usepackage[colorlinks=true,linkcolor=blue,citecolor=blue,urlcolor=blue]{hyperref}

\DeclareUnicodeCharacter{03A9}{\ensuremath{\Omega}}          
\DeclareUnicodeCharacter{211D}{\ensuremath{\mathbb{R}}}      
\DeclareUnicodeCharacter{2192}{\ensuremath{\rightarrow}}     
\DeclareUnicodeCharacter{2200}{\ensuremath{\forall}}         
\DeclareUnicodeCharacter{2208}{\ensuremath{\in}}             
\DeclareUnicodeCharacter{2209}{\ensuremath{\notin}}          
\DeclareUnicodeCharacter{2264}{\ensuremath{\leq}}            
\DeclareUnicodeCharacter{03B4}{\ensuremath{\delta}}          
\DeclareUnicodeCharacter{03B5}{\ensuremath{\varepsilon}}     
\DeclareUnicodeCharacter{03B7}{\ensuremath{\eta}}            
\DeclareUnicodeCharacter{03BC}{\ensuremath{\mu}}             
\DeclareUnicodeCharacter{03C1}{\ensuremath{\rho}}            
\DeclareUnicodeCharacter{03C6}{\ensuremath{\varphi}}         
\DeclareUnicodeCharacter{2115}{\ensuremath{\mathbb{N}}}      
\DeclareUnicodeCharacter{2203}{\ensuremath{\exists}}         
\DeclareUnicodeCharacter{2227}{\ensuremath{\wedge}}          
\DeclareUnicodeCharacter{2228}{\ensuremath{\vee}}            
\DeclareUnicodeCharacter{2260}{\ensuremath{\neq}}            
\DeclareUnicodeCharacter{2265}{\ensuremath{\geq}}            
\DeclareUnicodeCharacter{2286}{\ensuremath{\subseteq}}       
\DeclareUnicodeCharacter{2016}{\ensuremath{\Vert}}           
\DeclareUnicodeCharacter{2190}{\ensuremath{\leftarrow}}      
\DeclareUnicodeCharacter{2191}{\ensuremath{\uparrow}}        
\DeclareUnicodeCharacter{2194}{\ensuremath{\leftrightarrow}} 
\DeclareUnicodeCharacter{21A6}{\ensuremath{\mapsto}}         
\DeclareUnicodeCharacter{2202}{\ensuremath{\partial}}        
\DeclareUnicodeCharacter{221E}{\ensuremath{\infty}}          
\DeclareUnicodeCharacter{222B}{\ensuremath{\int}}            
\DeclareUnicodeCharacter{2A0D}{\ensuremath{\fint}}           
\DeclareUnicodeCharacter{1D4DD}{\ensuremath{\mathcal{N}}}    
\DeclareUnicodeCharacter{207B}{\ensuremath{^{-}}}            
\DeclareUnicodeCharacter{00B9}{\ensuremath{^{1}}}            
\DeclareUnicodeCharacter{1DA0}{\textsuperscript{f}}          
\DeclareUnicodeCharacter{2080}{\textsubscript{0}}            
\DeclareUnicodeCharacter{2081}{\textsubscript{1}}            
\DeclareUnicodeCharacter{2082}{\textsubscript{2}}            

\fvset{fontsize=\footnotesize, xleftmargin=10pt}
\DefineVerbatimEnvironment{lean}{Verbatim}{}

\newcommand{\R}{\mathbb{R}}
\newcommand{\eps}{\varepsilon}
\newcommand{\pnlap}{\Delta_p^{\mathrm{N}}}       
\newcommand{\inflap}{\Delta_\infty^{\mathrm{N}}} 
\newcommand{\Op}{T_{\eps}}                       
\newcommand{\Fop}{\mathcal{F}}                   

\DeclarePairedDelimiter{\abs}{\lvert}{\rvert}

\DeclarePairedDelimiterX{\inner}[2]{\langle}{\rangle}{#1, #2}

\makeatletter
\NewDocumentCommand{\fint}{}{\mathop{\mathpalette\fint@bar\relax}\!\int}
\newcommand*{\fint@bar}[2]{
  \sbox\z@{$\m@th#1\intop$}
  \rlap{\kern.225\wd\z@
        $\m@th#1\vcenter{\hrule\@height\fint@rule#1\@width.55\wd\z@}$}}
\newcommand*{\fint@rule}[1]{
  \ifx#1\scriptscriptstyle\fontdimen8\scriptscriptfont\thr@@
  \else\ifx#1\scriptstyle\fontdimen8\scriptfont\thr@@
  \else\fontdimen8\textfont\thr@@\fi\fi}
\makeatother

\DeclareUrlCommand\leanref{}

\theoremstyle{plain}
\newtheorem{theorem}{Theorem}[section]

\newtheorem*{maintheorem}{Main Theorem}
\theoremstyle{definition}
\newtheorem{definition}[theorem]{Definition}
\theoremstyle{remark}
\newtheorem{remark}[theorem]{Remark}

\begin{document}

\title[A formally verified convergence theorem for the tug-of-war scheme]
      {A formally verified convergence theorem\\
       for the tug-of-war scheme of the normalized $p$-Laplacian}

\author{Juan J. Manfredi}
\address{Department of Mathematics, University of Pittsburgh,
  Pittsburgh, PA 15260, USA}
\email{manfredi@pitt.edu}

\author{Behrooz Moosavi Ramezanzadeh}
\address{Department of Mathematics, University of Pittsburgh,
  Pittsburgh, PA 15260, USA}
\email{behroozmoosavi@pitt.edu}


\begin{abstract}
We present a complete formalization, using Lean~4, of the convergence of
$p$-harmonious functions to the viscosity solution of the Dirichlet problem
for the normalized $p$-Laplacian, $2 \le p < \infty$, as the step size goes to
zero. By $p$-harmonious functions we mean solutions of the dynamic programming
principle of tug-of-war with noise. Our proof follows the argument of Barles
and Souganidis and is formalized in full, depending only on Mathlib and the
standard axioms of Lean. In particular, we include the asymptotic expansion of
the scheme, the viscosity properties of the half-relaxed limits, existence and
uniqueness for the discrete problem, the comparison principle and uniqueness for the
homogeneous equation,
Aleksandrov's theorem, and Jensen's lemma.
\end{abstract}

\maketitle

\section{Introduction}\label{sec:intro}

For $1 < p < \infty$ the normalized $p$-Laplacian is
\begin{equation}\label{eq:pnlap}
  \pnlap u = \Delta u + (p-2)\,\inflap u, \qquad
  \inflap u = \Bigl\langle D^2u\,\tfrac{\nabla u}{\abs{\nabla u}},
                \tfrac{\nabla u}{\abs{\nabla u}}\Bigr\rangle .
\end{equation}
Here $u$ is a function defined on an open subset of $\R^d$, $d \ge 1$, with
gradient $\nabla u$ and Hessian $D^2u$. The Euclidean inner product is denoted
by $\langle\cdot,\cdot\rangle$, and \eqref{eq:pnlap} is defined where
$\nabla u \ne 0$.
The operator is homogeneous of degree one and not in divergence form, so the
natural notion of solution is the viscosity one. Its connection with
stochastic games goes back to Peres, Schramm, Sheffield and Wilson, who
identified the infinity Laplacian as the limit of tug-of-war \cite{PSSW09},
and to Peres and Sheffield, who reached finite $p$ by adding noise
\cite{PS08}. Manfredi, Parviainen and Rossi characterized $p$-harmonic
functions by an asymptotic mean value property \cite{MPR10} and introduced
$p$-harmonious functions, the solutions of the discrete dynamic programming
principle
\begin{equation}\label{eq:dpp}
  u_\eps(x) = \frac{\alpha}{2}\Bigl(\sup_{B_\eps(x)} u_\eps
     + \inf_{B_\eps(x)} u_\eps\Bigr)
     + \beta \fint_{B_\eps(x)} u_\eps, \qquad
  \alpha = \frac{p-2}{p+d}, \quad \beta = \frac{d+2}{p+d}.
\end{equation}
Here $\eps > 0$ is the step size, and $B_\eps(x)$ is the open ball of radius
$\eps$ centred at $x$. For a ball $B$ we denote by $\abs B$ its volume, and by
\[
  \fint_B v = \frac{1}{\abs B}\int_B v(y)\,dy
\]
the average of $v$ over $B$ with respect to Lebesgue measure. For $p \ge 2$
the weights satisfy $\alpha \ge 0$, $\beta > 0$ and $\alpha + \beta = 1$. In
the game, $\alpha$ is the probability of a tug-of-war move, in which the
choices of the two players produce the supremum and the infimum, and $\beta$
is the probability of a random move, uniformly distributed in $B_\eps(x)$,
which produces the average. Manfredi, Parviainen and Rossi showed that bounded
$p$-harmonious functions in bounded Lipschitz domains, with boundary values
prescribed by a fixed bounded continuous function $F$, converge as $\eps$ goes
to zero to the viscosity solution of the Dirichlet problem with boundary
values $F$ \cite{MPR12}.

\subsection*{Main result}
We report a complete formalization of a version of this convergence theorem,
for domains satisfying the exterior sphere condition, in Lean~4
\cite{dMU21} on top of the Mathlib library \cite{Mathlib20}.
The development builds with Lean~4.34.0-rc2 (toolchain
\texttt{leanprover/lean4:v4.34.0-rc2}) against Mathlib at commit
\texttt{85e3a25}.\footnote{The full commit hash is
\texttt{85e3a25e006c35636f0e53b0e9296caca2685bc0}.}

Recall that an open set $\Omega \subset \R^d$ satisfies the \emph{exterior
sphere condition} if for every $y \in \partial\Omega$ there is a ball
$B(z,\delta)$ with $B(z,\delta) \cap \Omega = \emptyset$ and
$y \in \partial B(z,\delta)$.

\begin{maintheorem}[\leanref{PLaplacian.pharmonious_exists_converge}]
Let $2 \le p < \infty$ and $d \ge 1$. Let $\Omega \subset \R^d$ be a bounded
open set satisfying the exterior sphere condition, and let
$F \colon \R^d \to \R$ be bounded and continuous. Then:
\begin{enumerate}
  \item for every $\eps > 0$ there is a unique function $u_\eps$, $p$-harmonious in
    $\Omega$ with data $F$ (Definition~\ref{def:pharm}), equal to $F$ outside
    $\Omega$, with
    $\sup\abs{u_\eps} \le \sup\abs F$;
  \item the half-relaxed limits $\overline u$ and $\underline u$ of this
    family (Section~\ref{ssec:relaxed}) coincide on $\overline\Omega$ and equal $F$ on $\partial\Omega$, and
    $u_\eps \to \overline u$ uniformly on $\overline\Omega$.
\end{enumerate}
The limit is a viscosity solution of $-\pnlap u = 0$ in $\Omega$
(Definition~\ref{def:visc}). The uniqueness of $u_\eps$ in part (1) is proved
in Theorem~\ref{thm:discunique}. Part (2)
holds for every uniformly bounded family of $p$-harmonious functions with data
$F$ that equal $F$ outside $\Omega$ (\leanref{pharmonious_converge}).
\end{maintheorem}

\noindent In Lean, the Main Theorem reads as follows; the exterior sphere
condition is the hypothesis \texttt{hext}, and the conclusion lists the
existence of the family, its bounds, and the three statements of part (2).
\begin{lean}
theorem pharmonious_exists_converge [Nontrivial (EuclSpace d)] (hd : 0 < d)
    {p : ℝ} (hp : 2 ≤ p) (hpd : 0 < p + (d : ℝ)) {Ω : Set (EuclSpace d)}
    (hΩo : IsOpen Ω) (hΩb : Bornology.IsBounded Ω) {F : EuclSpace d → ℝ}
    (hFc : Continuous F) {C : ℝ} (hFC : ∀ y, |F y| ≤ C)
    (hext : ∀ y ∈ frontier Ω, ∃ (z : EuclSpace d) (δ : ℝ), 0 < δ ∧
      ‖y - z‖ = δ ∧ ∀ x ∈ Ω, δ ≤ ‖x - z‖) :
    ∃ u : ℝ → EuclSpace d → ℝ,
      (∀ ε > 0, IsPHarmonious p ε Ω F (u ε)) ∧
      (∀ ε > 0, ∀ x, x ∉ Ω → u ε x = F x) ∧
      (∀ ε y, |u ε y| ≤ C) ∧
      (∀ x ∈ closure Ω, ubar u x = lbar u x) ∧
      (∀ x ∈ frontier Ω, ubar u x = F x) ∧
      (∀ η > 0, ∃ ε₀ > 0, ∀ ε, 0 < ε → ε < ε₀ → ∀ x ∈ closure Ω,
        |u ε x - ubar u x| ≤ η)
\end{lean}

The development contains no unproved statements. Lean reports that
\leanref{pharmonious_exists_converge} depends only on the axioms
\leanref{propext}, \leanref{Classical.choice} and \leanref{Quot.sound},
that is, on classical mathematics. The source is
available at
\url{https://github.com/BehroozMoosavi/tugofwar_formalization}, tag
\texttt{v1.3}; a technical companion \cite{MMcomp} documents every layer of
the development file by file.

\subsection*{Remarks}
We make the following remarks.

\emph{Scope.} The mathematical results are established; this paper describes
their formalization in Lean.

\emph{An analytic proof.} The dynamic programming principle is the definition
of the discrete problem;
existence follows the analytic argument of Luiro, Parviainen and Saksman
\cite{LPS14}, applied to the operator restricted to the domain, and
convergence follows as in del Teso, Manfredi and Parviainen \cite{dTMP22}.

\emph{Aleksandrov's theorem and Jensen's lemma.} Neither was available in
Mathlib. Our proof of Aleksandrov's theorem \cite{Ale39} obtains subgradients
from the proximal map, rather than from a supporting-hyperplane argument based
on the Hahn--Banach theorem, and derives the second-order expansion from the
subgradient inequalities by Riemann sums, rather than by integrating the
gradient along segments. Our proof of Jensen's lemma \cite{Jen88} uses neither
mollification nor the area formula.

\emph{Precise statements.} Several statements of the informal theory
had to be made precise; Section~\ref{sec:findings} collects a few examples.

\subsection*{Related work}
Formalized analysis of partial differential equations is relatively recent.
Van Doorn and Macbeth formalized the Gagliardo--Nirenberg--Sobolev inequality
\cite{vDM24}, and Armstrong and Kempe the interior De Giorgi--Nash--Moser
theory \cite{AK26}. Both concern divergence-form equations and weak solutions.
In probability, Marion formalized the Ionescu--Tulcea theorem and infinite
products of probability measures in Mathlib \cite{Mar26}. Manfredi has
also formalized in Lean the regularity of $p$-harmonic functions in the plane
of Iwaniec and Manfredi \cite{IMlean}, with Shirsho Mukherjee the comparison
principles for degenerate sub-elliptic equations in non-divergence form of
Manfredi and Mukherjee \cite{MMulean}, and interior $C^{1,\alpha}$ estimates
for infinity-harmonic functions \cite{InfHlean}.

\subsection*{Organization}
In Section~\ref{sec:setting} we fix the definitions, and in
Section~\ref{sec:plan} we give the plan of the proof. In
Sections~\ref{sec:consistency}--\ref{sec:boundary} we carry out its main
steps. Section~\ref{sec:findings} collects a few examples from the
formalization.

\section{Setting}\label{sec:setting}

Throughout the paper we fix the exponent $p$ with $2 \le p < \infty$, the
dimension $d \ge 1$, and the step size $\eps > 0$.

\subsection{The scheme}
A function $v \colon \R^d \to \R$ is \emph{admissible} if it is Borel
measurable and bounded. For admissible $v$ the step operator is
\begin{equation}\label{eq:op}
  \Op v(x) = \frac{\alpha}{2}\Bigl(\sup_{B_\eps(x)} v + \inf_{B_\eps(x)} v\Bigr)
    + \beta \fint_{B_\eps(x)} v ,
\end{equation}
with $\alpha,\beta$ as in \eqref{eq:dpp}. The balls are open
(see Section~\ref{sec:findings}). The operator is monotone, commutes with the
addition of constants, and is therefore nonexpansive in the supremum norm.
Monotonicity uses $p \ge 2$: for $p < 2$ the weight $\alpha$ is negative.

\begin{definition}\label{def:pharm}
Let $\Omega \subset \R^d$ and let $F \colon \R^d \to \R$ be bounded; $F$ is
the boundary data. An admissible $u$ is \emph{$p$-harmonious in $\Omega$ with
data $F$} if $u = \Op u$ on $\Omega$ and $u = F$ on the strip
$\Gamma_\eps = \{x \notin \Omega : \operatorname{dist}(x,\Omega) \le \eps\}$.
\end{definition}

\subsection{Viscosity solutions}
We write the equation as $\Fop(\nabla u, D^2u) = 0$, where
\[
  \Fop(q, X) = -\operatorname{tr} X - (p-2)\,
    \Bigl\langle X\tfrac{q}{\abs q}, \tfrac{q}{\abs q}\Bigr\rangle ,
  \qquad q \ne 0 .
\]
Here $q \in \R^d$, $X$ ranges over the symmetric $d \times d$ matrices, and
$\operatorname{tr} X$ is the trace of $X$. We also write
$\pnlap(q, X) = -\Fop(q, X)$, so that $\pnlap u = \pnlap(\nabla u, D^2u)$.
The operator $\Fop$ is degenerate elliptic, that is,
$\Fop(q, X) \le \Fop(q, Y)$ whenever $X - Y$ is positive semidefinite.

The function $\Fop$ is not defined at $q = 0$. To state the equation at points
where the gradient vanishes, we use the lower semicontinuous envelope
$\Fop_*$ of $\Fop$, the largest lower semicontinuous function below $\Fop$,
and the upper semicontinuous envelope $\Fop^*$, the smallest upper
semicontinuous function above $\Fop$. Both agree with $\Fop$ where $q \ne 0$.
For $p \ge 2$ their values at $q = 0$ are
\begin{equation}\label{eq:env}
  \Fop_*(0,X) = -\operatorname{tr}X - (p-2)\lambda_{\max}(X), \qquad
  \Fop^*(0,X) = -\operatorname{tr}X - (p-2)\lambda_{\min}(X),
\end{equation}
where $\lambda_{\min}(X)$ and $\lambda_{\max}(X)$ denote the smallest and the
largest eigenvalue of $X$. We use these formulas directly, so that no general
envelope construction is needed.

\begin{definition}\label{def:visc}
An upper semicontinuous $u$ is a \emph{viscosity subsolution} of
$-\pnlap u = 0$ in an open set $\Omega$ if, whenever $\varphi \in C^2$ and
$u - \varphi$ has a local maximum at $x_0 \in \Omega$,
\[
  \Fop(\nabla\varphi(x_0), D^2\varphi(x_0)) \le 0 \ \text{ if } \nabla\varphi(x_0)
  \ne 0, \qquad \Fop_*(0, D^2\varphi(x_0)) \le 0 \ \text{ otherwise}.
\]
A lower semicontinuous $u$ is a \emph{viscosity supersolution} of
$-\pnlap u = 0$ in $\Omega$ if, whenever $\varphi \in C^2$ and $u - \varphi$
has a local minimum at $x_0 \in \Omega$,
\[
  \Fop(\nabla\varphi(x_0), D^2\varphi(x_0)) \ge 0 \ \text{ if } \nabla\varphi(x_0)
  \ne 0, \qquad \Fop^*(0, D^2\varphi(x_0)) \ge 0 \ \text{ otherwise}.
\]
A function is a \emph{viscosity solution} if it is both a viscosity
subsolution and a viscosity supersolution; in particular, it is continuous.
\end{definition}

\subsection{Jets and second-order expansions}
Let $u$ be a real function defined near $x_0$. We say that $u$ has a
\emph{second-order expansion} $(q, X)$ at $x_0$, with $q \in \R^d$ and
$X$ a symmetric matrix, if
\[
  u(x_0 + h) = u(x_0) + \langle q, h\rangle + \tfrac12\langle Xh, h\rangle
  + o(\abs h^2) \qquad (h \to 0).
\]
Here $o(\abs h^2)$ denotes a remainder $r(h)$ with $r(h)/\abs h^2 \to 0$ as
$h \to 0$; for this notion, see \cite[Section~6.4]{EG15}. If $u \in C^2$, this is the second-order Taylor expansion, with
$q = \nabla u(x_0)$ and $X = D^2u(x_0)$; the definition, however, does not
require $u$ to be differentiable near $x_0$. The pair $(q, X)$ is a
\emph{superjet} of $u$ at $x_0$ if
\[
  u(x_0+h) \le u(x_0)+\langle q,h\rangle
    +\tfrac12\langle Xh,h\rangle+o(\abs h^2)
  \qquad (h\to0).
\]
The pair is a \emph{subjet} if the reverse inequality holds.
Definition~\ref{def:visc} can equivalently
be stated with jets: its conditions then apply to superjets (for
subsolutions) and subjets (for supersolutions) $(q, X)$ in place of
$(\nabla\varphi(x_0), D^2\varphi(x_0))$.

\subsection{Half-relaxed limits}\label{ssec:relaxed}
For a uniformly bounded family $(u_\eps)$, we define its upper and lower
half-relaxed limits by
\[
  \overline u(x) = \limsup_{\eps \to 0^+,\ y \to x} u_\eps(y), \qquad
  \underline u(x) = \liminf_{\eps \to 0^+,\ y \to x} u_\eps(y).
\]
The limits are taken along the product of the punctured right neighbourhood
of $0$ with the neighbourhood of $x$. They are defined on all of $\R^d$ and
finite, $\underline u \le \overline u$, and $\overline u$ is upper and
$\underline u$ lower semicontinuous. Since $-u_\eps$ is $p$-harmonious with
data $-F$ whenever $u_\eps$ is (the operator is odd),
$\underline u = -\overline{(-u)}$, and every statement about $\underline u$
below is read off from the corresponding statement about $\overline u$.

\section{Plan of the proof}\label{sec:plan}

Let $(u_\eps)_{\eps > 0}$ be the family given by part (1) of the Main Theorem:
for each $\eps > 0$, the function $u_\eps$ is $p$-harmonious in $\Omega$ with
data $F$ and equals $F$ outside $\Omega$. Let $\overline u$ and
$\underline u$ be its half-relaxed limits. The argument of Barles and
Souganidis \cite{BS91} reduces part (2) of the Main Theorem to three statements:
\begin{enumerate}
  \item[(A)] The upper half-relaxed limit $\overline u$ is a viscosity
    subsolution, and the lower half-relaxed limit $\underline u$ is a
    viscosity supersolution in $\Omega$ (Section~\ref{sec:relaxed}).
  \item[(B)] The boundary inequalities $\overline u \le F \le \underline u$
    hold on $\partial\Omega$
    (Section~\ref{sec:boundary}).
  \item[(C)] The comparison principle holds: a bounded upper semicontinuous
    subsolution lies below a bounded lower semicontinuous supersolution in
    $\Omega$ if it does so on $\partial\Omega$ (Section~\ref{sec:cp}).
\end{enumerate}
Part (1) of the Main Theorem, the existence and uniqueness of $u_\eps$ for each
$\eps > 0$, is proved separately: existence by iteration, and uniqueness by a
propagation argument (Section~\ref{ssec:existence}). Given (A), (B) and (C),
$\overline u \le \underline u$ in $\Omega$, hence equality on
$\overline\Omega$; the common value is continuous, and uniform convergence
follows by a compactness argument against the semicontinuity of the relaxed
limits. Statement (A) rests on the asymptotic expansion of the scheme
(Section~\ref{sec:consistency}).

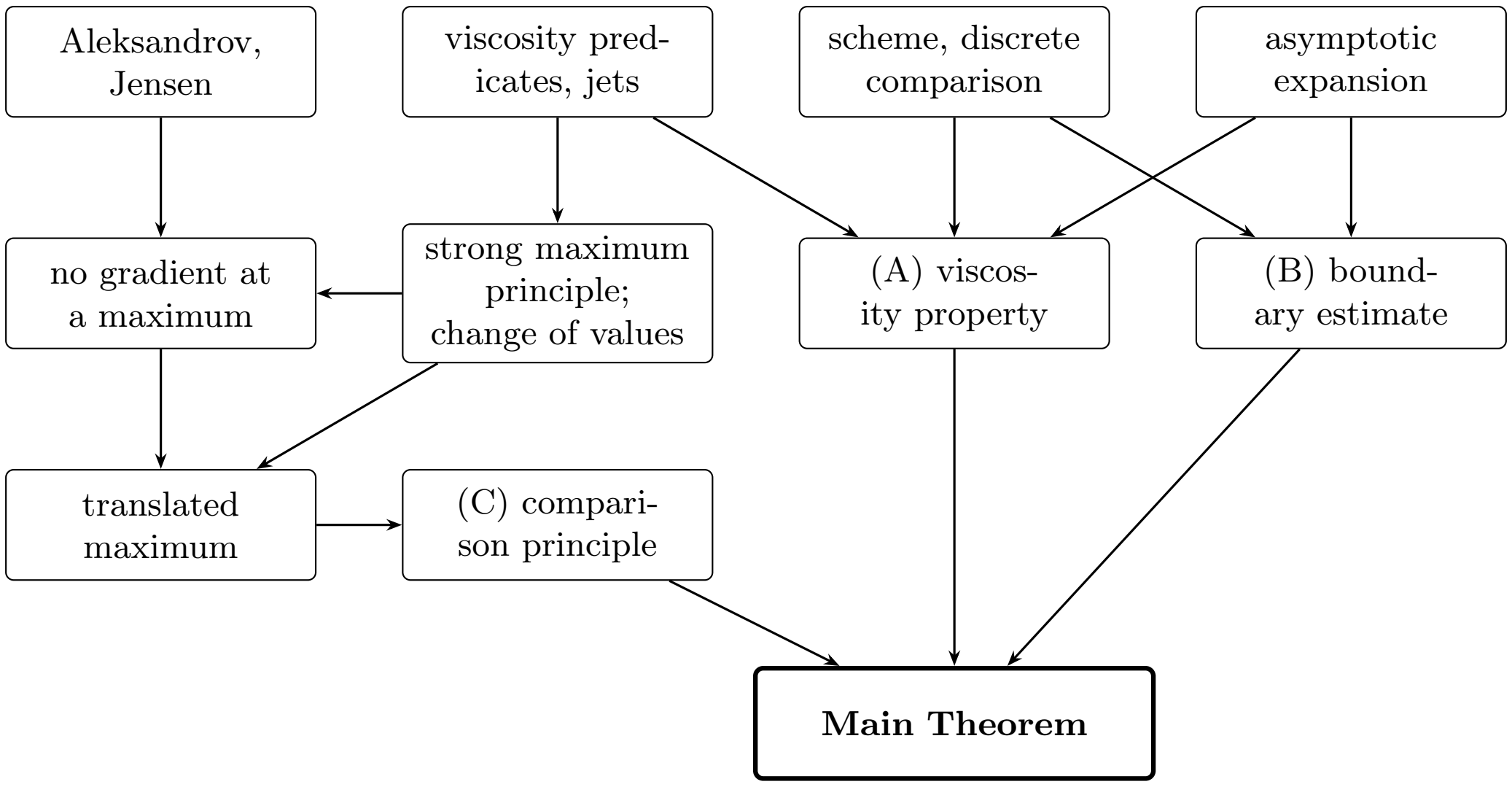
\begin{figure}[ht]
\centering
\begin{tikzpicture}[
  box/.style={draw, rounded corners=2pt, align=center, font=\footnotesize,
              text width=2.6cm, minimum height=1.0cm, inner sep=3pt},
  main/.style={box, very thick, text width=3.4cm},
  arr/.style={-{Stealth[length=4pt]}, semithick}]
  \node[box] (AJ) at (0,0)    {Aleksandrov, Jensen};
  \node[box] (L2) at (3.6,0)  {viscosity predicates, jets};
  \node[box] (L1) at (7.2,0)  {scheme, discrete comparison};
  \node[box] (L3) at (10.8,0) {asymptotic expansion};
  \node[box] (S2) at (0,-2.1)   {no gradient at a maximum};
  \node[box] (SC) at (3.6,-2.1) {strong maximum principle; change of values};
  \node[box] (V4) at (7.2,-2.1) {(A) viscosity property};
  \node[box] (BD) at (10.8,-2.1) {(B) boundary estimate};
  \node[box] (S3) at (0,-4.2)   {translated maximum};
  \node[box] (CP) at (3.6,-4.2) {(C) comparison principle};
  \node[main] (MT) at (7.2,-6.0) {\textbf{Main Theorem}};
  \draw[arr] (AJ) -- (S2);  \draw[arr] (L2) -- (SC);
  \draw[arr] (L2) -- (V4);  \draw[arr] (L1) -- (V4);
  \draw[arr] (L3) -- (V4);  \draw[arr] (L1) -- (BD);
  \draw[arr] (L3) -- (BD);  \draw[arr] (SC) -- (S2);
  \draw[arr] (S2) -- (S3);  \draw[arr] (SC) -- (S3);
  \draw[arr] (S3) -- (CP);  \draw[arr] (CP) -- (MT);
  \draw[arr] (V4) -- (MT);  \draw[arr] (BD) -- (MT);
\end{tikzpicture}
\caption{Dependencies among the main results; an arrow means that the target
uses the source.}\label{fig:chain}
\end{figure}

The formal development is organized in four layers, summarized in
Table~\ref{tab:layers}; Figure~\ref{fig:chain} shows how the main results
depend on each other.

\begin{table}[ht]
\centering\small
\begin{tabular}{@{}llr@{}}
\toprule
Layer & Content & Lines (approx.) \\
\midrule
1 & the scheme: operator, measurability, existence, uniqueness, comparison
  & 3\,300 \\
2 & viscosity framework: generator, predicates, consistency, stability
  & 1\,100 \\
3 & the asymptotic expansion of $\Op$, with an explicit error rate
  & 5\,500 \\
4 & relaxed limits, comparison principle, boundary, convergence
  & 19\,100 \\
\midrule
 & total (225 files) & 29\,400 \\
\bottomrule
\end{tabular}
\caption{The four layers of the development. Layer~4 includes
Aleksandrov's theorem and Jensen's lemma.}\label{tab:layers}
\end{table}

\section{Consistency: the asymptotic expansion}\label{sec:consistency}

\begin{theorem}[\leanref{expansion_uniform}, \leanref{eventually_bound_le}]\label{thm:exp}
Let $\varphi \in C^3$ with $\nabla\varphi(x_0) \ne 0$. There is $\rho > 0$
such that, for every $y \in B_\rho(x_0)$ and every $0 < \eps \le \rho$,
\[
  \Op\varphi(y) = \varphi(y) + \frac{\eps^2}{2(p+d)}\,\pnlap\varphi(y)
     + E_\eps(y),
  \qquad\text{where}\qquad
  \sup_{y \in B_\rho(x_0)} \abs{E_\eps(y)} = o(\eps^2) \quad (\eps \to 0^+).
\]
\end{theorem}

\noindent In Lean, the first declaration gives the expansion with an explicit
bound on the error, uniform in $y$, and the second shows that this bound is
eventually at most $\eta\eps^2$ for every $\eta > 0$.
\begin{lean}
theorem expansion_uniform [Nontrivial (EuclSpace d)] (hd : 0 < d)
    {p : ℝ} (hp : 2 ≤ p) (hpd : 0 < p + (d : ℝ))
    {φ : EuclSpace d → ℝ} (hφ : ContDiff ℝ 3 φ) {x₀ : EuclSpace d}
    (hx₀ : gradient φ x₀ ≠ 0) :
    ∃ ρ > 0, ∃ H K Ka g : ℝ, 0 < g ∧ 0 ≤ K ∧ 0 ≤ H ∧ 0 ≤ Ka ∧
      ∀ y ∈ ball x₀ ρ, ∀ ε : ℝ, 0 < ε → ε ≤ ρ →
        ∃ E : ℝ,
          tugOfWarOp p ε φ y
              = φ y + ε ^ 2 / (2 * (p + (d : ℝ)))
                  * normPLapOf p (gradient φ y)
                      (fderiv ℝ (fderiv ℝ φ) y) + E
            ∧ |E| ≤ (H * ε * Real.sqrt
                  (2 * ε * ((H * ε ^ 2 + 2 * K * ε ^ 3) / g))
                + K * ε ^ 3) + Ka * ε ^ 3

theorem eventually_bound_le {H K Ka g η : ℝ} (hη : 0 < η) :
    ∀ᶠ ε in nhdsWithin (0 : ℝ) (Set.Ioi 0),
      (H * ε * Real.sqrt (2 * ε * ((H * ε ^ 2 + 2 * K * ε ^ 3) / g))
          + K * ε ^ 3) + Ka * ε ^ 3 ≤ η * ε ^ 2
\end{lean}

The right-hand side of \eqref{eq:op} is the sum of two terms: the mean of the
supremum and the infimum of $v$ over $B_\eps(x)$, multiplied by $\alpha$, and
the average of $v$ over $B_\eps(x)$, multiplied by $\beta$. We expand each
term for $v = \varphi$.

The average yields the Laplacian through the second moments of Lebesgue
measure on the ball,
\[
  \fint_{B_\eps(0)} z_i z_j \, dz = \frac{\eps^2}{d+2}\,\delta_{ij},
\]
where $\delta_{ij}$ is the Kronecker delta. This identity is exact.
Averaging the second-order Taylor expansion gives
\[
  \fint_{B_\eps(x)} \varphi = \varphi(x) + \frac{\eps^2}{2(d+2)}\,
  \Delta\varphi(x) + O(\eps^3).
\]
The first moments vanish, and the remainder of the second-order Taylor
expansion of $\varphi \in C^3$ is bounded by $C\eps^3$ on $B_\eps(x)$.
The constant $C$ is independent of $\eps$ and can be chosen uniformly for
$x$ in a fixed sufficiently small neighbourhood of $x_0$. Thus the
$O(\eps^3)$ remainder is also $o(\eps^2)$, since $C\eps^3/\eps^2\to0$.
The mean of the supremum and the infimum yields the infinity
Laplacian: the maximizer of $\varphi$ on $\overline B_\eps(x)$ lies within
$O(\eps^{3/2})$ of $x + \eps\,\nabla\varphi(x)/\abs{\nabla\varphi(x)}$, the
first-order terms of the two extreme values cancel, and
\[
  \tfrac12\Bigl(\sup_{B_\eps(x)}\varphi + \inf_{B_\eps(x)}\varphi\Bigr)
  = \varphi(x) + \frac{\eps^2}{2}\,\inflap\varphi(x) + o(\eps^2).
\]
The weights then combine as
\[
  \alpha\,\frac{\eps^2}{2}\,\inflap\varphi
  + \beta\,\frac{\eps^2}{2(d+2)}\,\Delta\varphi
  = \frac{\eps^2}{2(p+d)}\,\pnlap\varphi,
\]
the factor $d+2$ cancelling against the numerator of $\beta$.

Theorem~\ref{thm:exp} is uniform in the base point $y$, as the viscosity
property requires, since near-maximizers move with $\eps$. The error bound is
explicit: the Taylor remainder is cubic in the increment, and the
$O(\eps^{3/2})$ distance of an extremizer from the corresponding point
$y\pm\eps\nabla\varphi(y)/\abs{\nabla\varphi(y)}$ contributes
$O(\eps^{5/2})$ to the function-value remainder \cite[\S5]{MMcomp}.
At critical points the
extreme-value term is controlled instead by the eigenvalues of the Hessian,
which is where the envelopes \eqref{eq:env} appear.

\section{The viscosity property of the relaxed limits}\label{sec:relaxed}

\begin{theorem}[\leanref{ubar_isViscositySubsolution}, \leanref{lbar_isViscositySupersolution}]\label{thm:A}
If $(u_\eps)$ is uniformly bounded and each $u_\eps$ is $p$-harmonious in an
open set $\Omega$, then $\overline u$ is a viscosity subsolution of
$-\pnlap u = 0$ in $\Omega$, in the sense of Definition~\ref{def:visc}, and
$\underline u$ is a viscosity supersolution.
\end{theorem}

\begin{samepage}
\noindent In Lean:
\begin{lean}
theorem ubar_isViscositySubsolution [Nontrivial (EuclSpace d)] (hd : 0 < d)
    {p : ℝ} (hp : 2 ≤ p) (hpd : 0 < p + (d : ℝ))
    {Ω : Set (EuclSpace d)} (hΩ : IsOpen Ω) {F : EuclSpace d → ℝ}
    {u : ℝ → EuclSpace d → ℝ} {C : ℝ} (hC : ∀ ε y, |u ε y| ≤ C)
    (hu : ∀ ε > 0, IsPHarmonious p ε Ω F (u ε)) :
    IsViscositySubsolution p Ω (fun _ => 0) (ubar u)

theorem lbar_isViscositySupersolution [Nontrivial (EuclSpace d)] (hd : 0 < d)
    {p : ℝ} (hp : 2 ≤ p) (hpd : 0 < p + (d : ℝ))
    {Ω : Set (EuclSpace d)} (hΩ : IsOpen Ω) {F : EuclSpace d → ℝ}
    {u : ℝ → EuclSpace d → ℝ} {C : ℝ} (hC : ∀ ε y, |u ε y| ≤ C)
    (hu : ∀ ε > 0, IsPHarmonious p ε Ω F (u ε)) :
    IsViscositySupersolution p Ω (fun _ => 0) (lbar u)
\end{lean}
\end{samepage}

The statement uses $C^2$ test functions, arbitrary local maxima, and both the
regular and the singular conditions. The proof has three main steps.

\emph{Near-maximizers without extraction.} Suppose $\overline u - \varphi$
has a strict maximum at $x_0$. For each $\eps$ choose $x_\eps$ nearly
maximizing $u_\eps - \varphi$ on a ball about $x_0$, measured against that
parameter's own supremum. Applying the monotone operator to
$u_\eps \le \varphi + c_\eps$ and using $u_\eps = \Op u_\eps$ gives
\[
  \Op\varphi(x_\eps)-\varphi(x_\eps)\ge-\eta_\eps,
  \qquad \eta_\eps=o(\eps^2).
\]
The uniform expansion turns this into the subsolution
inequality in the limit. That the near-maximizers approach $x_0$ is shown by
a compactness argument on an annulus, using strictness, rather than by
extracting a convergent sequence.

\emph{The singular case without a gradient.} When $\nabla\varphi(x_0) = 0$,
the consistency estimate is unavailable. Instead, the extreme values of a
quadratic over a ball are compared with its extreme eigenvalue by reflecting
the ball through its centre, which gives the inequality involving
$\lambda_{\max}$ in \eqref{eq:env} directly.

\emph{From smooth to $C^2$ test functions.} The core argument is carried out
for $C^3$ test functions at strict extrema. The passage to $C^2$ functions and
non-strict extrema adds a quartic $\abs{x - x_0}^4$ and replaces $\varphi$ by
its second-order Taylor polynomial. This needs neither mollification nor the
symmetry of the Hessian, because every operator in the predicates sees the
Hessian only through its quadratic form $X(v,v)$.

\section{The comparison principle}\label{sec:cp}

\begin{theorem}[\leanref{visc_comparison}]\label{thm:C}
Let $\Omega \subset \R^d$ be bounded and open, and let $U, V \colon \R^d \to \R$
be bounded, with $U$ upper semicontinuous and $V$ lower semicontinuous on
$\R^d$. Suppose that $U$ is a viscosity subsolution and $V$ a viscosity
supersolution of $-\pnlap u = 0$ in $\Omega$. If $U(y) \le V(y)$ for every
$y \in \partial\Omega$, then $U \le V$ in $\Omega$.
\end{theorem}

\noindent In Lean, where the boundedness hypotheses appear as bounds on the
ranges:
\begin{lean}
theorem visc_comparison (hd : 0 < d) {p : ℝ} (hp : 2 ≤ p) {Ω : Set (EuclSpace d)}
    (hΩo : IsOpen Ω) (hΩb : Bornology.IsBounded Ω) {U V : EuclSpace d → ℝ}
    (hU : UpperSemicontinuous U) (hV : LowerSemicontinuous V)
    (hUb : BddAbove (range U)) (hUb' : BddBelow (range U))
    (hVb : BddBelow (range V)) (hVb' : BddAbove (range V))
    (hsub : IsViscositySubsolution p Ω (fun _ => 0) U)
    (hsup : IsViscositySupersolution p Ω (fun _ => 0) V)
    (hbd : ∀ x ∈ frontier Ω, U x ≤ V x) :
    ∀ x ∈ Ω, U x ≤ V x
\end{lean}

The boundary condition is pointwise. Since $U$ and $V$ are semicontinuous on
all of $\R^d$, it implies the usual boundary inequality
\[
  \limsup_{x \to y,\ x \in \Omega} U(x) \le U(y) \le V(y)
  \le \liminf_{x \to y,\ x \in \Omega} V(x) \qquad (y \in \partial\Omega).
\]
In the application $U = \overline u$ and $V = \underline u$, which are defined
and semicontinuous on all of $\R^d$.

Applying Theorem~\ref{thm:C} in both directions gives uniqueness.

\begin{theorem}[\leanref{visc_unique}]\label{thm:unique}
Two bounded continuous viscosity solutions of $-\pnlap u = 0$ in a bounded open
set $\Omega$ that agree on $\partial\Omega$ agree in $\Omega$. In particular,
the limit in the Main Theorem does not depend on the family $(u_\eps)$.
\end{theorem}

\noindent In Lean:
\begin{lean}
theorem visc_unique (hd : 0 < d) {p : ℝ} (hp : 2 ≤ p)
    {Ω : Set (EuclSpace d)} (hΩo : IsOpen Ω) (hΩb : Bornology.IsBounded Ω)
    {u₁ u₂ : EuclSpace d → ℝ}
    (h₁c : Continuous u₁) (h₂c : Continuous u₂)
    (h₁a : BddAbove (range u₁)) (h₁b : BddBelow (range u₁))
    (h₂a : BddAbove (range u₂)) (h₂b : BddBelow (range u₂))
    (h₁sub : IsViscositySubsolution p Ω (fun _ => 0) u₁)
    (h₁sup : IsViscositySupersolution p Ω (fun _ => 0) u₁)
    (h₂sub : IsViscositySubsolution p Ω (fun _ => 0) u₂)
    (h₂sup : IsViscositySupersolution p Ω (fun _ => 0) u₂)
    (hbd : ∀ x ∈ frontier Ω, u₁ x = u₂ x) :
    ∀ x ∈ Ω, u₁ x = u₂ x
\end{lean}

\subsection{The obstruction}
For a strict subsolution, $-\pnlap U \le -c < 0$, the classical doubling
argument with the theorem of sums \cite{CIL92} applies, and we formalize it
first. The relaxed limits, however, satisfy the equation with $c = 0$.
Adding $\frac\delta2\abs x^2$ to the subsolution shifts the Hessian by
$\delta I$ but also the gradient by $\delta x$, and in the doubling argument
the resulting error in the \emph{direction} of the gradient is of order
$\delta\abs x/\abs{x - y}$, unbounded as the doubled points merge. A monotone
change of unknown keeps the direction, but the strictness it creates vanishes
at critical points. Another route to comparison uses the equivalence of
viscosity solutions with weak solutions \cite{JLM01}; the argument below uses
only viscosity theory.

\subsection{The change of values}
Following Barles and Busca \cite{BB01}, let $\psi(s) = s - \eta s^2$ with
$\eta > 0$ small, and $W = \psi\circ V$.
If $V$ has the second-order expansion $(q, X)$ at a point, as at the
points given by Aleksandrov's theorem \cite[Section~6.4]{EG15}, then $W$
has the second-order expansion
\[
  (cq,\ cX-2\eta\,q\otimes q), \qquad c=\psi'(V)>0,
\]
where
$q \otimes q$ is the matrix $(q_i q_j)_{i,j}$, and the
supersolution inequality gives
\begin{equation}\label{eq:cv}
  \pnlap\bigl(cq,\ cX - 2\eta\, q\otimes q\bigr) \le -2\eta(p-1)\abs q^2 .
\end{equation}
This is strict exactly where $q \ne 0$. The composition is used only at
points where $V$ is twice differentiable, so no viscosity chain rule is
needed.

\subsection{No maximum with nonzero gradient}
After regularization (Section~\ref{ssec:reg}), $U$ is semiconvex and $V$ and
$W$ are semiconcave: there are constants $\lambda, \mu \ge 0$ such that
$U + \frac\lambda2\abs{\cdot}^2$ is convex and $W - \frac\mu2\abs{\cdot}^2$ is
concave. A vector $q$ is a \emph{lower support} of $U$ at $x$ if
\[
  U(x + h) \ge U(x) + \langle q, h\rangle - \tfrac\lambda2\abs h^2
  \qquad\text{for all } h;
\]
by semiconvexity, $U$ has a lower support at every point. Let $x$ be a
local maximum of $U - W$. We show that every lower support $q$ of $U$ at $x$
vanishes. Suppose not. Jensen's lemma, applied
to $U - W - \delta\abs{\cdot - x}^2$, gives a set of positive measure of points
at which small linear perturbations attain their maximum; by Aleksandrov's
theorem $U$ and $-V$ are twice differentiable almost everywhere, so the set
contains a point $y_0$ of twice differentiability of both. There the
first- and second-order conditions at a maximum, together with
\eqref{eq:cv} and the subsolution inequality for $U$, give
\[
  0 \le -\kappa + \tfrac{\kappa}{3} + \tfrac{\kappa}{3},
  \qquad \kappa = 2\eta(p-1)\Bigl(\frac{\abs q}{2C}\Bigr)^2 > 0,
\]
where $C$ is an upper bound for $\psi' \circ V$; this is a contradiction. A
quantitative gradient estimate, $\abs{a - q} \le (\lambda + 2\mu)\abs{y - x}$
for lower supports $a$ at nearby points $y$, allows all parameters to be fixed
in advance, so that a single point of the positive-measure set suffices and
no sequence of points is extracted.

\subsection{The translated maximum}
On a connected open set $G$ with $U - W$ larger inside than on $\partial G$,
consider $M(h) = \max_{\overline G}\bigl(U(\cdot + h) - W\bigr)$. For small $h$
its maximizers are interior, and the previous step applies to the translate:
its lower support vanishes. Hence $M(h + k) \ge M(h) - \frac\lambda2\abs k^2$,
and subdividing a segment into $N$ steps loses at most
$\frac{\lambda}{2N}\abs k^2$, so $M$ is constant near $0$. On the contact set
$K = \{U - W = M(0)\}$, then, $U$ has local maxima and $V$ local minima; the
strong maximum and minimum principles make both locally constant, so $K$ is
open. It is also relatively closed and nonempty, so $K = G$, which contradicts
the gap at the boundary.

The strong maximum principle is proved with the barrier
$(\abs{x - y_0}^2 - B)^2$ on a thin annulus, and is stated for functions whose
superjets with nonzero gradient satisfy the subsolution inequality. In that
form it applies to the regularizations directly, without proving that they
are viscosity subsolutions.

\subsection{Regularization}\label{ssec:reg}
We replace $U$ by its sup-convolution and $V$ by its inf-convolution, which
are semiconvex and semiconcave; by the magic property of sup-convolutions
\cite[Lemma~A.5]{CIL92} their jets satisfy the jet inequalities wherever the
maximizers of the convolution lie in $\Omega$. This needs $U$ and $V$ bounded
on both sides. Semicontinuity and compactness give a band near
$\partial\Omega$ on which the regularized difference is small, and the
connected component of a superlevel set containing a point where
$U - V$ is positive supplies the set $G$ of the previous step.

\subsection{Aleksandrov's theorem and Jensen's lemma}
Aleksandrov's theorem gives the following second-order expansion
\cite[Section~6.4]{EG15}.
\begin{theorem}[\leanref{alexandrov}]\label{thm:alex}
A convex function $f$ on $\R^d$ has a second-order expansion at almost every
point: for almost every $x$ there are $q \in \R^d$ and a
symmetric bilinear form $X$ with
\[
  f(x + h) = f(x) + \langle q, h\rangle + \tfrac12 X(h,h) + o(\abs h^2)
  \qquad (h \to 0).
\]
\end{theorem}

\noindent In Lean, where \texttt{AlexandrovThm d} is the proposition stated
above:
\begin{lean}
theorem alexandrov : AlexandrovThm d
\end{lean}

For convex $f$, the proximal map $J(y) = \operatorname{argmin}_x\bigl(f(x) +
\frac12\abs{x - y}^2\bigr)$ is one-Lipschitz and onto, and $y - J(y)$ is a
subgradient at $J(y)$; subgradients at every point are obtained as limits of
rescaled proximal differences, without a separation theorem. Rademacher's
theorem and the area formula, both in Mathlib, show that almost every $x$ is
of the form $J(y)$ with $DJ(y)$ invertible. At such points the subgradient
has a linear approximation, and the second-order expansion follows from the
subgradient inequalities alone, by a Riemann sum for $\int_0^1 t\,dt$; no
integral is used.

\begin{theorem}[\leanref{jensen_positive_measure}]\label{thm:jensen}
Let $r>0$ and $\delta>0$, and let $W$ be semiconvex on $\R^n$.
Suppose that $W$ is larger at some point
of $\overline B(x_0, r)$ than $2\delta r$ plus its values on the sphere. Then
the set of points of $B(x_0,r)$ at which some $W + \langle \xi, \cdot\rangle$
with $\xi \in \R^n$, $\abs \xi \le \delta$, attains its maximum over the
closed ball $\overline B(x_0,r)$ has positive measure.
\end{theorem}

\begin{samepage}
\noindent In Lean:
\begin{lean}
theorem jensen_positive_measure {n : ℕ} (W : EuclSpace n → ℝ) (C : ℝ)
    (hC : 0 ≤ C) (hconv : ConvexOn ℝ Set.univ (fun x => W x + C / 2 * ‖x‖ ^ 2))
    (x₀ xs : EuclSpace n) (r δ : ℝ) (hr : 0 < r) (hδ : 0 < δ)
    (hxs : xs ∈ closedBall x₀ r)
    (hsphere : ∀ y, ‖y - x₀‖ = r → W y + 2 * δ * r < W xs) :
    0 < volume {x | x ∈ ball x₀ r ∧ ∃ p : EuclSpace n, ‖p‖ ≤ δ ∧
      ∀ y ∈ closedBall x₀ r, W y + inner ℝ p y ≤ W x + inner ℝ p x}
\end{lean}
\end{samepage}

The tilt at such a point is determined by a subgradient map $g$, which sends
the set onto a ball and is Lipschitz on the parts of the set away from the
sphere, where $W$ is squeezed between its semiconvex lower bound and the upper
bound given by the maximum. A null set would have a null image.

\section{Boundary behaviour, existence and convergence}\label{sec:boundary}

\begin{theorem}[\leanref{ubar_le_bdry}]\label{thm:B}
Under the hypotheses of the Main Theorem, $\overline u(y_0) \le F(y_0)$ at
every $y_0 \in \partial\Omega$; by negation, $F(y_0) \le \underline u(y_0)$.
\end{theorem}

\noindent In Lean, where the exterior ball at $y_0$ is $B(z_0, \delta_0)$:
\begin{lean}
theorem ubar_le_bdry [Nontrivial (EuclSpace d)] (hd : 0 < d) {p : ℝ} (hp : 2 ≤ p)
    (hpd : 0 < p + (d : ℝ)) {Ω : Set (EuclSpace d)} {F : EuclSpace d → ℝ}
    {u : ℝ → EuclSpace d → ℝ} {C : ℝ} (hC : ∀ ε y, |u ε y| ≤ C)
    (hu : ∀ ε > 0, IsPHarmonious p ε Ω F (u ε))
    (hout : ∀ ε > 0, ∀ x, x ∉ Ω → u ε x = F x)
    {y₀ : EuclSpace d} (hFc : ContinuousAt F y₀)
    {z₀ : EuclSpace d} {δ₀ : ℝ} (hδ₀ : 0 < δ₀) (hy₀z : ‖y₀ - z₀‖ = δ₀)
    (hext : ∀ x ∈ Ω, δ₀ ≤ ‖x - z₀‖) :
    ubar u y₀ ≤ F y₀
\end{lean}

The proof compares $u_\eps$ with a barrier through a discrete comparison
principle: \emph{if $u = \Op u$ and $\Op v \le v - \kappa$ on a set $D$, for
some $\kappa > 0$, and $u \le v$ within $\eps$ outside $D$, then $u \le v$ on
$D$.} Unlike comparison of two solutions, this needs no attained maximum.

\emph{The exterior ball.} Let $B(z_0,\delta_0)$ be an exterior ball at $y_0$,
so that $\abs{y_0 - z_0} = \delta_0$ and $\abs{x - z_0} \ge \delta_0$ for every
$x \in \overline\Omega$. We replace it by the ball of half the radius whose
centre is the midpoint of $y_0$ and $z_0$,
\[
  z = \tfrac12(y_0 + z_0), \qquad \delta = \tfrac12\delta_0 .
\]
Apollonius's identity expresses the distance to a midpoint through the
distances to the two endpoints: for every $x \in \R^d$,
\[
  \abs{x - z}^2 = \tfrac12\abs{x - y_0}^2 + \tfrac12\abs{x - z_0}^2
    - \tfrac14\abs{y_0 - z_0}^2 ,
\]
which follows by expanding the squares (it is the parallelogram law). For
$x \in \overline\Omega$ the middle term is at least $\frac12\delta_0^2$ and the
last equals $\frac14\delta_0^2$, so
\begin{equation}\label{eq:apollonius}
  \abs{x - z}^2 \ge \delta^2 + \tfrac12\abs{x - y_0}^2
  \qquad (x \in \overline\Omega).
\end{equation}
Thus $y_0$ is the only point of $\overline\Omega$ on the sphere
$\partial B(z,\delta)$, and the points of $\overline\Omega$ with
$\abs{x - z}^2 \le \delta^2 + \tau$ lie within $\sqrt{2\tau}$ of $y_0$: the
comparison region shrinks to $y_0$ as $\tau \to 0$. For the original ball,
every point of the sphere $\partial B(z_0,\delta_0)$ could lie on
$\partial\Omega$, and no such localization would hold.

\emph{The barrier.} With $t = \abs{x - z}^2 - \delta^2$, let
\[
  \psi = t - Mt^2, \qquad M = \frac{2(d + p - 2) + 1}{8(p-1)\,\delta^2}.
\]
On the shell $0 \le t \le 1/(4M)$ the gradient of $\psi$ does not vanish and
$\pnlap\psi \le -1$, so by the expansion, uniform on the compact shell,
\[
  \Op(c + A\psi) \le c + A\psi - \frac{A\eps^2}{4(p+d)}
\]
on the shell, for all constants $c \in \R$ and $A > 0$ and all sufficiently
small $\eps$. The barrier is a polynomial, hence globally smooth as the
expansion requires; the classical barrier $\abs{x - z}^{-\sigma}$ is not. By
\eqref{eq:apollonius}, $t \ge \frac12\abs{x - y_0}^2 \ge 0$ on
$\overline\Omega$, and the comparison is made on the part of $\Omega$ near
$y_0$ where $t$ is small.

\smallskip
\noindent\emph{Note.} Under the exterior sphere condition used here, this
proof is simpler than the corresponding argument of del Teso, Manfredi
and Parviainen \cite{dTMP22}.

\subsection{Existence with exterior data}\label{ssec:existence}
\begin{theorem}[\leanref{exists_pharmonious}]\label{thm:exist}
Let $\Omega$ be bounded and measurable, $\eps > 0$, and $F$ measurable with
$\abs F \le C$. Then there is a $p$-harmonious function in $\Omega$ with data
$F$, equal to $F$ outside $\Omega$ and bounded by $C$.
\end{theorem}

\noindent In Lean:
\begin{lean}
theorem exists_pharmonious (p ε : ℝ) (hp : 2 ≤ p) (hε : 0 < ε) (hd : 0 < d)
    (hpd : 0 < p + (d : ℝ)) {Ω : Set (EuclSpace d)} (hΩm : MeasurableSet Ω)
    (hΩb : Bornology.IsBounded Ω)
    {F : EuclSpace d → ℝ} (hFm : Measurable F) {C : ℝ}
    (hFC : ∀ y, |F y| ≤ C) :
    ∃ u : EuclSpace d → ℝ, IsPHarmonious p ε Ω F u
      ∧ (∀ x, x ∉ Ω → u x = F x) ∧ ∀ y, |u y| ≤ C
\end{lean}

The proof iterates the operator restricted to $\Omega$, $\widetilde T u = \Op u$
on $\Omega$ and $\widetilde T u = u$ elsewhere, starting from $u_0 = -C$ on
$\Omega$ and $u_0 = F$ outside. The iterates increase, stay bounded by $C$, and
are never moved outside $\Omega$. For uniform convergence we use the argument of
\cite{LPS14}: with $M$ the limit of $\sup_{\Omega}(u - u_j)$, $u$ the limit and
$x_0 \in \Omega$ a point nearly attaining the supremum,
\[
  M - 2\delta \le \Op u_\ell(x_0) - \Op u_k(x_0)
  \le \alpha \sup_{B_\eps(x_0)}(u - u_k) + \beta \fint_{B_\eps(x_0)}(u - u_k)
  \le \alpha(M + \delta) + \delta
\]
for suitable indices $k < \ell$, the average being small by dominated
convergence and the independence of the ball volume from its centre. Hence
$\beta M \le (3 + \alpha)\delta$ for every $\delta > 0$, and $M = 0$. Because the
deficit vanishes outside $\Omega$, the convergence is uniform on all of $\R^d$,
and the limit is a fixed point of $\widetilde T$ by nonexpansiveness.

\begin{theorem}[\leanref{pharmonious_unique}]\label{thm:discunique}
Let $\Omega$ be bounded. Two $p$-harmonious functions in $\Omega$ with data $F$
that both equal $F$ outside $\Omega$ are equal.
\end{theorem}

\noindent In Lean:
\begin{lean}
theorem pharmonious_unique [Nontrivial (EuclSpace d)] (p ε : ℝ) (hp : 2 ≤ p)
    (hε : 0 < ε) (hpd : 0 < p + (d : ℝ)) {Ω : Set (EuclSpace d)}
    (hΩb : Bornology.IsBounded Ω) {F u v : EuclSpace d → ℝ}
    (hu : IsPHarmonious p ε Ω F u) (hv : IsPHarmonious p ε Ω F v)
    (huF : ∀ x, x ∉ Ω → u x = F x) (hvF : ∀ x, x ∉ Ω → v x = F x) :
    u = v
\end{lean}

We follow the strategy of Luiro, Parviainen and Saksman
\cite[Theorem~2.2]{LPS14}: near-maxima of the difference of two solutions
propagate in a fixed direction until they leave the bounded domain. Let
$w = v - u$ and $M = \sup w$, and suppose $M > 0$. Since $w = 0$ outside
$\Omega$, every point where $w > 0$ lies in $\Omega$. If $w(x) > M - \delta$,
then the dynamic programming principle gives
\[
  w(x) \le \alpha M + \beta \fint_{B_\eps(x)} w,
  \qquad\text{hence}\qquad
  \fint_{B_\eps(x)} w > M - \frac{\delta}{\beta}.
\]
Comparing integrals of $M - w$ over $B_\eps(x)$ and over the ball
$S = B(x + \frac34\eps e, \frac14\eps)$, for a fixed unit vector $e$, shows
that $w > M - C\delta$ at some point of $S$, which lies at least $\eps/2$
further in the direction $e$; the constant $C$ depends only on $d$, $p$ and
$\eps$. After $N$ such steps, with $N\eps/2$ larger than the width of
$\Omega$, the point has left $\Omega$, yet $w > M - C^N\delta_0 > 0$ there if
$\delta_0 = M/(2C^N)$: a contradiction.

\smallskip
\noindent\emph{Note.} In \cite{LPS14} the maximum is first
shown to be attained almost everywhere on a ball, by compactness, and the
exact maximal set is propagated; our proof propagates near-maxima with
geometrically growing errors instead, which avoids the compactness step. The
resulting argument is simpler than the one in \cite{LPS14}.

With (A), (B), (C) and Theorems~\ref{thm:exist} and~\ref{thm:discunique} in
hand, the Main Theorem follows as described in Section~\ref{sec:plan}.

\begin{remark}
The Main Theorem assumes $u_\eps = F$ outside $\Omega$, not only on the strip
$\Gamma_\eps$. This is needed: the relaxed limits at a boundary point see
values at points farther than $\eps$ from $\Omega$, where the dynamic
programming principle imposes nothing.
\end{remark}

\section{Examples from the formalization}\label{sec:findings}

This section collects a few examples of statements and estimates that had
to be made precise in the formalization. The companion \cite[\S9]{MMcomp}
gives further details.

\subsection{Open balls and measurability}
If the extreme values in \eqref{eq:op} are taken over \emph{closed} balls,
the operator need not preserve Borel measurability: there are bounded Borel
functions whose supremum over closed balls is not Borel. The admissible class
is then not invariant under the operator.
Over open balls, the set
\[
  \{x : \sup_{B_\eps(x)} v > \lambda\}
\]
is a union of open balls, which gives measurability.

\subsection{The restricted operator}
Iterating $\Op$ on all of $\R^d$ gives a fixed-point equation $\Op u=u$
everywhere, without imposing the boundary data. The estimate of
\cite{LPS14} controls the deficit only near
$\Omega$, while passing the operator through the limit needs uniform
convergence on every ball used by the operator. With the operator restricted
to $\Omega$, the deficit vanishes
outside $\Omega$ and the same estimate closes the argument.

\subsection{The boundary strip}
Every exterior point sampled by an open ball of radius $\eps$ centred in
$\Omega$ lies at distance strictly less than $\eps$ from $\Omega$. The definition above
prescribes data on the strip with distance at most $\eps$; prescribing
values also at distance exactly $\eps$ does not affect the operator
inside $\Omega$.

\subsection{Analytic existence and uniqueness arguments}
The proofs of \cite{MPR12} use the tug-of-war game, with supermartingales and
optional stopping. We follow the analytic arguments of \cite{LPS14}. In the
existence argument, the uniformity needed follows from the independence of the
volume of a ball from its centre.

\subsection{Two representations of higher derivatives}
Mathlib represents second derivatives in two ways: as the curried derivative
of the derivative and as the flat iterated derivative. The Taylor estimates
of Section~\ref{sec:consistency} use explicit identities between these
representations for differentiation and norm estimates.

\subsection{Polynomial barriers}
To avoid computing the Hessians of explicit functions, both barriers of the
proof enter only through exact polynomial expansions: for the strong maximum
principle, the expansion turns a maximum into a superjet; for the boundary
barrier, uniqueness of the coefficients in a second-order expansion identifies the operator value with
its value at the explicit expansion.

\subsection{Quantifiers in error estimates}
For a fixed $\eps>0$, the inequality
\[
  \abs E\le K\eps^2
\]
holds for any $E$ with $K=\abs E/\eps^2$. An error estimate for the expansion
in Theorem~\ref{thm:exp} requires constants independent of $\eps$, so the
bounds must be established uniformly in that parameter. For its
$o(\eps^2)$ conclusion, the error bound must be at most $\eta\eps^2$ for
every $\eta>0$, once $\eps$ is sufficiently small, uniformly in the base
point. In the development,
the propositions used by subsequent results were stated before their proofs.
For example, Aleksandrov's theorem and Jensen's lemma were hypotheses of the
comparison principle until both were proved.

\subsection{How the development was produced}
As in \cite{AK26}, the formalization was carried out with extensive use of a
large language model, which proposed formal statements and drafted Lean files
under the authors' direction. Every draft was compiled against the pinned
toolchain, and corrections were driven by the compiler's output. The
mathematical route was fixed by the authors. The Lean code is checked by
Lean's kernel, and the Main Theorem depends only on the standard axioms.

\section{Conclusion}\label{sec:conclusion}

We have formalized existence, uniqueness and convergence for the tug-of-war
scheme of the normalized $p$-Laplacian, including the comparison principle for
the limit equation, in about 29\,400 lines of Lean.

\section*{Acknowledgements}
We thank the Mathlib community, on whose library this development is built,
and the developers of Lean. J.~J.~Manfredi was partially supported by Simons
Foundation gift no.~962828.

\bibliographystyle{amsplain}
\bibliography{ref}

\end{document}